\documentclass[12pt]{article}
\usepackage{amssymb, amsfonts, amsmath, amsthm}

\def\3{\subset }
\def\4{\subseteq }
\def\ov{\overline}

\def\cala{{\cal A}}

\def\0{\leqno}
\def\a{{\alpha}}
\def\barr{\begin{array}}
\def\earr{\end{array}}
\def\dd{\displaystyle}

\def\Z{{\rlap{$\kern2pt{\rm Z}$}{\rm Z}\,}}

\title{\bf Solitary quotients of finite groups}
\author{Marius T\u arn\u auceanu}
\date{October 1, 2012}

\begin{document}

\maketitle

\begin{abstract}
In this paper we introduce and study the lattice of normal
subgroups of a group $G$ that determine solitary quotients. It is
closely connected to the well-known lattice of solitary subgroups
of $G$ (see \cite{5}). A precise description of this lattice is
given for some particular classes of finite groups.
\end{abstract}



\section{Introduction}

The relation between the structure of a group and the structure of
its lattice of subgroups constitutes an important domain of
research in group theory. The topic has enjoyed a rapid
development starting with the first half of the �20 century. Many
classes of groups determined by different properties of partially
ordered subsets of their subgroups (especially lattices of
subgroups) have been identified. We refer to Suzuki's book
\cite{9}, Schmidt's book \cite{8} or the more recent book
\cite{11} by the author for more information about this theory.

The starting point for our discussion is given by the paper
\cite{5}, where the lattice Sol($G$) of solitary subgroups of a
group $G$ has been introduced. A natural idea is to study the
normal subgroups of $G$ that induce solitary quotients. The set of
these subgroups also forms a lattice, denoted by QSol($G$), which
constitutes a "dual" for Sol($G$). The first steps in studying
this new lattice associated to a group is the main goal of the
current paper.

In the following, given a group $G$, we will denote by $L(G)$ the
subgroup lattice of $G$. Recall that $L(G)$ is a complete bounded
lattice with respect to set inclusion, having initial element the
trivial subgroup $\{1\}$ and final element $G$, and its binary
operations $\wedge, \vee$ are defined by
$$H\wedge K=H\cap K,\ H\vee K=\langle H\cup K\rangle, \mbox{ for all } H,K\in
L(G).$$Two important modular sublattices of $L(G)$ are the normal
subgroup lattice and the characteristic subgroup lattice of $G$,
usually denoted by $N(G)$ and ${\rm Char}(G)$, respectively. Note
also that ${\rm Char}(G)$ is a sublattice of $N(G)$.

Two groups $G_1$ and $G_2$ will be called \textit{L-isomorphic} if
their subgroup lattices $L(G_1)$ and $L(G_2)$ are isomorphic. A
\textit{duality} from $G_1$ onto $G_2$ is a bijective map $\delta:
L(G_1)\longrightarrow L(G_2)$ such that the following equivalent
conditions are satisfied:
\begin{itemize}
\item[--] $\forall\hspace{1mm} H, K\in L(G_1)$, we have $H\leq K$ if and only if $K^{\delta}\leq H^{\delta}$;
\item[--] $(H\wedge K)^{\delta}=H^{\delta}\vee K^{\delta}$, for all $H, K\in L(G_1)$;
\item[--] $(H\vee K)^{\delta}=H^{\delta}\wedge K^{\delta}$, for all $H, K\in L(G_1)$.
\end{itemize}We say that a group $G$
\textit{has a dual} if there exists a duality from $G$ to some
group $\ov{G}$, and that $G$ is \textit{self-dual} if there exists
a duality from $G$ onto $G$.

A subgroup $H\in L(G)$ is called a \textit{solitary subgroup} of
$G$ if $G$ does not contain another isomorphic copy of $H$, that
is
$$\forall\hspace{1mm} K\in L(G), K\cong H \Longrightarrow
K=H.$$By Theorem 25 of \cite{5}, the set Sol($G$) consisting of
all solitary subgroups of $G$ forms a lattice with respect to set
inclusion, which is called \textit{the lattice of solitary
subgroups of} $G$.

The paper is organized as follows. In Section 2 we present some
basic properties and results on the lattice QSol($G$) associated
to a finite group $G$. Section 3 deals with this lattice for
finite abelian groups. The most important theorem in this section
gives a complete description of the lattices Sol($G$) and
QSol($G$) for such a group $G$. In the final section some
conclusions and further research directions are indicated.

Most of our notation is standard and will usually not be repeated
here. Basic definitions and results on lattices and groups can be
found in \cite{1,2} and \cite{3,4,10}, respectively. For subgroup
lattice concepts we refer the reader to \cite{8,9,11}.

\section{The lattice QSol($G$)}

Let $G$ be a finite group and QSol($G$) be the set of normal
subgroups of $G$ that determine solitary quotients, that is
$${\rm QSol}(G)=\{H\in N(G) \mid \forall\hspace{1mm} K\in N(G), G/K\cong G/H \Longrightarrow
K=H\}.$$We remark that all elements of QSol($G$) are
characteristic subgroups of $G$. It is also obvious that $\{1\}$
and $G$ are contained in QSol($G$). If QSol($G$) consists only of
these two subgroups, then we will call $G$ a \textit{quotient
solitary free group}.

\bigskip\noindent{\bf Remark.} Clearly, another important
subgroup of $G$ that belongs to QSol($G$) is $G'$. We infer that
if $G$ is quotient solitary free, then $G'=G$ or $G'=\{1\}$ and
therefore $G$ is either perfect or abelian. Observe also that the
elementary abelian $p$-groups are quotient solitary free.
\bigskip

Our first result shows that QSol($G$) can naturally be endowed
with a lattice structure.

\bigskip\noindent{\bf Proposition 2.1.} {\it Let $G$ be a finite group.
Then {\rm QSol($G$)} is a lattice with respect to set inclusion.}
\bigskip

\noindent{\bf Proof.} Let $H_1, H_2\in{\rm QSol}(G)$ and $H\in
N(G)$ such that $G/H\cong G/H_1\cap H_2$. Choose an isomorphism
$f: G/H_1\cap H_2\longrightarrow G/H$. Then there are two normal
subgroups $H_1', H_2'$ of $G$ with $f(H_i/H_1\cap H_2)=H_i'/H$,
$i=1,2$. It is easy to see that $H_1'\cap H_2'=H$. One obtains
$$\dd\frac{G}{H_i}\hspace{1mm}\cong\hspace{1mm}\dd\frac{\dd\frac{G}{H_1\cap H_2}}{\dd\frac{H_i}{H_1\cap H_2}}\hspace{1mm}\cong\hspace{1mm}\dd\frac{\dd\frac{G}{H}}{\dd\frac{H_i'}{H}}\hspace{1mm}\cong\hspace{1mm}\dd\frac{G}{H_i'}\hspace{1mm},$$which
implies that $H_i'=H_i$, $i=1,2$, and so $H=H_1\cap H_2$.
Therefore $H_1\cap H_2$ is the meet of $H_1$ and $H_2$ in
QSol($G$).

Obviously, the join of $H_1$ and $H_2$ in QSol($G$) also exists,
and consequently QSol($G$) is a lattice. Note that $\{1\}$ and $G$
are, respectively, the initial element and the final element of
this lattice.
\hfill\rule{1,5mm}{1,5mm}
\bigskip

An exhaustive description of the above lattice for an important
class of finite groups, the dihedral groups
$$D_{2n}=\langle x,y\mid x^n=y^2=1,\ yxy=x^{-1}\rangle, \hspace{1mm}n\geq 3,$$is indicated
in the following.

\bigskip\noindent{\bf Example.} The structure of $L(D_{2n})$ is
well-known: for every divisor $r$ of $n$, $D_{2n}$ possesses a
unique cyclic subgroup of order $r$ (namely $\langle
x^{\frac{n}{r}}\rangle $) and $\frac{n}{r}$ subgroups isomorphic
to $D_{2\frac{n}{r}}$ (namely $\langle x^{\frac{n}{r}},x^iy\rangle
$, \,\,\,$i=0,1,...,\frac{n}{r}-1$). Remark that $D_{2n}$ has
always a maximal cyclic normal subgroup $M=\langle x\rangle
\cong\Z_n$. Clearly, all subgroups of $M$ are normal in $D_{2n}$.
On the other hand, if $n$ is even, then $D_{2n}$ has another two
maximal normal subgroups of order $n$, namely $M_1=\langle
x^2,y\rangle $ and $M_2=\langle x^2,xy\rangle $, both isomorphic
to $D_n$. In this way, one obtains
$$N(D_{2n})=\left\{\barr{lll}
L(M) \cup \{D_{2n}\},&n\equiv1\hspace{1mm}({\rm mod} \hspace{1mm}2)\\
\\
L(M) \cup \{D_{2n}, M_1, M_2\},&n\equiv0 \hspace{1mm} ({\rm mod}
\hspace{1mm}2).\earr\right.$$We easily infer that
$${\rm QSol}(D_{2n})=\left\{\barr{lll}
L(M) \cup \{D_{2n}\},&n\equiv1\hspace{1mm}({\rm mod} \hspace{1mm}2)\\
\\
L(M)^* \cup \{D_{2n}\},&n\equiv0 \hspace{1mm} ({\rm mod}
\hspace{1mm}2),\earr\right.$$where $L(M)^*$ denotes the set of all
proper subgroups of $M$.
\bigskip

Observe that the equality between the lattices Sol($G$) and
QSol($G$) associated to a finite group $G$ fails. For example, we
have $\langle x^2\rangle \in$ {\rm QSol($D_8$)} but $\langle
x^2\rangle \notin$ {\rm Sol($D_8$)}, respectively $\langle
x\rangle \in$ {\rm Sol($D_8$)} but $\langle x\rangle \notin$ {\rm
QSol($D_8$)}.
\bigskip

By looking to the dihedral groups $D_{2n}$ with $n$ odd, we also
infer that there exist finite groups $G$ such that QSol($G$)=
$N(G)$. Other examples of such groups are the finite groups
without normal subgroups of the same order and, in particular, the
finite groups $G$ for which $N(G)$ is a chain (as simple groups,
symmetric groups, cyclic $p$-groups or finite groups of order $p^n
q^m$ ($p,q$ distinct primes) with cyclic Sylow subgroups and
trivial center -- see Exercise 3, page 497, \cite{8}). Note also
that for a finite group $G$ the condition QSol($G$)= $N(G)$ is
very close to the conditions of Theorem 9.1.6 of \cite{8} that
characterize the distributivity of $N(G)$.
\bigskip

The following proposition shows that the relation "to be quotient
solitary" is transitive.

\bigskip\noindent{\bf Proposition 2.2.} {\it Let $G$ be a finite group and $H\supseteq K$
be two normal subgroups of $G$. If $H\in$ {\rm QSol($G$)} and
$K\in$ {\rm QSol($H$)}, then $K\in$ {\rm QSol($G$)}.}
\bigskip

\noindent{\bf Proof.} Let $K_1\in N(G)$ such that $G/K_1\cong G/K$
and take an isomorphism $f: G/K\longrightarrow G/K_1$. Set
$H_1/K_1=f(H/K)$, where $H_1$ is a normal subgroup of $G$
containing $K_1$. It follows that
$$\dd\frac{G}{H}\hspace{1mm}\cong\hspace{1mm}\dd\frac{\dd\frac{G}{K}}{\dd\frac{H}{K}}\hspace{1mm}\cong\hspace{1mm}\dd\frac{\dd\frac{G}{K_1}}{\dd\frac{H_1}{K_1}}\hspace{1mm}\cong\hspace{1mm}\dd\frac{G}{H_1}\hspace{1mm},$$which
leads to $H_1=H$. One obtains $H/K_1\cong H/K$ and therefore
$K_1=K$. Hence $K\in$ {\rm QSol($G$)}. \hfill\rule{1,5mm}{1,5mm}
\bigskip

Next we will study the connections between the lattices {\rm
QSol($G$)} and {\rm QSol($\ov{G}$)}, where $\ov{G}$ is an
epimorphic image of $G$. We mention that a proper subgroup $H\in$
{\rm QSol($G$) will be called \textit{maximal} in {\rm QSol($G$)
if it is not properly contained in any proper subgroup of {\rm
QSol($G$).

\bigskip\noindent{\bf Proposition 2.3.} {\it Let $G$ be a finite group and $H$ be a
proper normal subgroup of $G$. Set $\ov{G}=G/H$ and denote by
$\pi: G\longrightarrow\ov{G}$ the canonical homomorphism. Then,
for every $K\in$ {\rm QSol($G$)}, we have $\pi(K)\in$ {\rm
QSol($\ov{G}$)}. In particular, if $H\in$ {\rm QSol($G$)} and
$\ov{G}$ is quotient solitary free, then $H$ is maximal in {\rm
QSol($G$)}.}
\bigskip

\noindent{\bf Proof.} Let $K\in$ {\rm QSol($G$)} and $K_1$ be a
normal subgroup of $G$ which contains $H$ and satisfies
$\ov{G}/\pi(K_1)\cong \ov{G}/\pi(K)$. Then $G/K_1\cong G/K$ and
thus $K_1=K$, in view of our hypothesis. Hence $\pi(K_1)=\pi(K)$.

Suppose now that $\ov{G}$ is quotient solitary free and let $K\in$
{\rm QSol($G$)} with $K\neq G$ and $H\subset K$. Then, by what we
proved above, $\ov{K}=\pi(K)\in$ {\rm QSol($\ov{G}$)} and
$\{\ov{1}\}\subset\ov{K}\subset\ov{G}$, a contradiction.
\hfill\rule{1,5mm}{1,5mm}

\bigskip\noindent{\bf Remark.} Under the hypotheses of Proposition 2.3,
$\pi$ fails to induce a bijection between the sets
$\cala=\{K\in{\rm QSol}(G)\mid H\subseteq K\}$ and QSol($\ov{G}$),
as follows by taking $G=D_{12}$ and $H$ the (unique) normal
subgroup of order 3 of $D_{12}$ (in this case $\cala$ consists of
three elements, namely $H$, $D_{12}$ and the cyclic subgroup of
order 6 in $D_{12}$, while
$\ov{D_{12}}\cong\mathbb{Z}_2\times\mathbb{Z}_2$ is quotient
solitary free).
\bigskip

In the following we assume that $G$ is nilpotent and let
$$G=\prod_{i=1}^k G_i\0(1)$$be the decomposition of $G$ as a direct
product of Sylow subgroups. Then the lattice $N(G)$ is
decomposable. More precisely every $H\in N(G)$ can be written as
$\prod_{i=1}^k H_i$ with $H_i\in N(G_i)$, $i=1,2,...,k$. We infer
that $H\in$ {\rm QSol($G$)} if and only if $H_i\in$ {\rm
QSol($G_i$)}, for all $i=\ov{1,k}$. In this way, the lattice {\rm
QSol($G$)} is also decomposable
$${\rm QSol}(G)\cong\prod_{i=1}^k {\rm QSol}(G_i)\0(2)$$and its study
is reduced to $p$-groups.
\bigskip

We first remark that for a finite $p$-group $G$, the lattice ${\rm
Char}(G)$ is in general strictly contained in {\rm QSol($G$)}
(take, for example, $G=S_{2^n}$, the quasi-dihedral group of order
$2^n$; being isomorphic to $\mathbb{Z}_{2^{n-1}}$, $D_{2^{n-1}}$
and $S_{2^{n-1}}$, respectively, the maximal subgroups of
$S_{2^n}$ are characteristic, but clearly they do not belong to
{\rm QSol($G$)}). In fact, a maximal subgroup $M$ of $G$ is
contained in {\rm QSol($G$)} if and only if $G$ is cyclic.
\bigskip

The following result will play an essential role in proving the
main theorem of Section 3. It illustrates another important
element of {\rm QSol($G$)} in the particular case of $p$-groups.

\bigskip\noindent{\bf Proposition 2.4.} {\it Let $G$ be a finite $p$-group.
Then $\Phi(G)$ is a maximal element of {\rm QSol($G$)}.}
\bigskip

\noindent{\bf Proof.} If $G/H\cong G/\Phi(G)$ for some $H\in
N(G)$, then $G/H$ is elementary abelian. Since $\Phi(G)$ is
minimal in $G$ with the property to determine an elementary
abelian quotient, we have $\Phi(G)\subseteq H$. On the other hand,
we know that $\Phi(G)$ and $H$ are of the same order. These lead
to $H=\Phi(G)$, that is $\Phi(G)\in {\rm QSol}(G)$.

We already have seen that $G/\Phi(G)$ is quotient solitary free.
According to Proposition 2.3, this implies the maximality of
$\Phi(G)$ in {\rm QSol($G$)}. \hfill\rule{1,5mm}{1,5mm}
\bigskip

Obviously, by Propositions 2.2 and 2.4 we infer that the Frattini
series of a finite $p$-group $G$ is contained in {\rm QSol($G$)},
that is
$$\{\Phi_n(G)\mid n\in\mathbb{N}\}\subseteq {\rm QSol}(G),$$where
$\Phi_0(G)=G$ and $\Phi_n(G)=\Phi(\Phi_{n-1}(G))$, for all $n\geq
1$. This can na\-tu\-ra\-lly be extended to finite nilpotent
groups.

\bigskip\noindent{\bf Corollary 2.5.} {\it Let $G$ be a finite nilpotent group. Then
$$\{\Phi_n(G)\mid n\in\mathbb{N}\}\subseteq {\rm QSol}(G).\0(3)$$Moreover, under the above notation,
the maximal elements of {\rm QSol($G$)} are
$$\Phi(G_i)\prod_{j\neq i} G_j\hspace{1mm}, \hspace{2mm}
i=1,2,...,k.$$}

In particular, we also obtain the following corollary.

\bigskip\noindent{\bf Corollary 2.6.} {\it Let $G$ be a finite nilpotent group. Then
{\rm QSol($G$)}= $N(G)$ if and only if $G$ is cyclic.}

\section{The case of finite abelian groups}

In this section we will focus on describing the lattice {\rm
QSol($G$)} associated to a finite abelian group $G$. As we have
seen above, it suffices to consider finite abelian $p$-groups. In
the following our main goal is to prove that for such a group the
relation (3) becomes an equality.
\bigskip

Recall first a famous theorem due to Baer (see, for example,
Theorem 8.1.4 of \cite{8} or Theorem 4.2 of \cite{9}) which states
that every finite abelian group (and, in particular, every finite
abelian $p$-group) $G$ is self-dual. Moreover, by fixing an
autoduality $\delta$ of $G$, we have
$$H\cong G/\delta(H)\mbox{ and } \delta(H)\cong G/H, \mbox{ for all } H\in
L(G).\0(4)$$These isomorphisms easily lead to the following
proposition.

\bigskip\noindent{\bf Proposition 3.1.} {\it Let $G$ be a finite abelian
$p$-group and $\delta$ be an autoduality of $G$. Then $\delta({\rm
QSol}(G))={\rm Sol}(G)$ and $\delta({\rm Sol}(G))={\rm QSol}(G)$,
that is $\delta$ induces an anti-isomorphism between the lattices
{\rm QSol($G$)} and {\rm Sol($G$)}. Moreover, we have
$\delta^2(H)=H$, for all $H\in {\rm QSol}(G)$.}
\bigskip

By Proposition 2.4, we know that $\Phi(G)$ is a maximal element of
QSol($G$). The $\Phi$-subgroup of a finite abelian $p$-group
satisfies some other simple but important properties in QSol($G$).

\bigskip\noindent{\bf Lemma 3.2.} {\it Let $G$ be a finite abelian
$p$-group and $G\cong\prod_{i=1}^k\Z_{p^{\a_i}}$ be the primary
decomposition of $G$. Then, for every proper subgroup $H$ of
\,{\rm QSol($G$)}, we have\,{\rm:}
\begin{itemize}
\item[{\rm a)}] $H\subseteq\Phi(G)${\rm;}
\item[{\rm b)}] $H\in {\rm QSol}(\Phi(G))$.
\end{itemize}}
\smallskip

\noindent{\bf Proof.} a) We shall proceed by induction on $k$. The
inclusion is trivial for $k=1$. Assume now that it holds for any
abelian $p$-group of rank less than $k$ and put $G=G_1\times G_2$,
where $G_1\cong\prod_{i=1}^{k-1}\Z_{p^{\a_i}}$ and
$G_2\cong\Z_{p^{\a_k}}$. According to Suzuki \cite{10}, vol. I,
(4.19), a subgroup $H$ of $G$ is uniquely determined by two
subgroups $H_1\subseteq H'_1$ of $G_1$, two subgroups
$H_2\subseteq H'_2$ of $G_2$ and an isomorphism $\varphi:
H'_1/H_1\longrightarrow H'_2/H_2$ (more exactly, we have
$H=\{(x_1,x_2)\in H'_1\times H'_2\mid\varphi(x_1H_1)=x_2H_2\})$.
Mention that $H'_i=\pi_i(H)$, $i=1,2$, where $\pi_1$ and $\pi_2$
are the projections of $G$ onto $G_1$ and $G_2$, respectively.
Clearly, Proposition 2.3 implies that $H'_i$ belongs to
QSol($G_i$), $i=1,2$. Since $H\in{\rm QSol}(G)$, it follows that
each $H'_i$ is properly contained in $G_i$. Indeed, if $H'_1=G_1$,
then $\pi_1$ induces a surjective homomorphism from $H$ to $G_1$,
and therefore $H$ has a quotient isomorphic to $G_1$. By duality,
it also has a subgroup isomorphic to $G_1$. This implies that
$G/H$ is isomorphic to a quotient of $G_2$, i.e. it is cyclic.
Thus $\delta(H)$ is a cyclic solitary subgroup of $G$. In other
words, $G$ contains a unique non-trivial cyclic subgroup of a
certain order, a contradiction. Similarly, we have $H'_2\neq G_2$.
Then, by the inductive hypothesis, one obtains
$H'_i\subseteq\Phi(G_i)$, $i=1,2$, and so
$$H\subseteq H'_1\times
H'_2\subseteq\Phi(G_1)\times\Phi(G_2)=\Phi(G).$$

b) By using Lemma 8.1.6 of \cite{8}, we infer that $\delta$
induces a bijection between the set of proper subgroups of
QSol($G$) contained in $\Phi(G)$ and ${\rm
Sol}(G/\delta(\Phi(G))$. Since the groups $G/\delta(\Phi(G))$ and
$\Phi(G)$ are isomorphic, their lattices of solitary subgroups are
also isomorphic. Finally, on account of Proposition 3.1, the
lattices ${\rm Sol}(\Phi(G))$ and ${\rm QSol}(\Phi(G))$ are
anti-isomorphic. Hence there is a bijection between the sets
$\{H\in{\rm QSol}(G)\mid H\subseteq \Phi(G)\}$ and ${\rm
QSol}(\Phi(G))$. On the other hand, by Proposition 2.2 we have
$${\rm QSol}(\Phi(G))\subseteq\{H\in{\rm QSol}(G)\mid H\subseteq
\Phi(G)\},$$and therefore these sets are equal. This completes the
proof.
\hfill\rule{1,5mm}{1,5mm}

\bigskip\noindent{\bf Remark.} An alternative way of proving a) of
Lemma 3.2 is obtained by using the lattice of characteristic
subgroups of $G$. According to Theorem 3.7 of \cite{7}, ${\rm
Char}(G)$ has a unique minimal element, say $M$, and clearly this
is solitary in $G$. It follows that $\delta(M)$ is the unique
maximal element of QSol($G$) and thus it will coincide with
$\Phi(G)$. In this way, all proper subgroups of QSol($G$) are
contained in $\Phi(G)$.
\bigskip

Let $r$ be the length of the Frattini series of $G$ (note that if
$G\cong\prod_{i=1}^k \Z_{p^{\a_i}}$,
$\a_1\leq\a_2\leq\cdots\leq\a_k$, is the primary decomposition of
$G$, then $r=\a_k$). By using Lemma 3.2 and a standard induction
on $r$, we easily come up with the conclusion that QSol($G$)
coincides with this series. The lattice Sol($G$) is also
completely determined in view of Proposition 3.1. Hence we have
proved the following theorem.

\bigskip\noindent{\bf Theorem 3.3.} {\it Let $G$ be a finite
abelian $p$-group and $G\cong\prod_{i=1}^k \Z_{p^{\a_i}}$ be the
primary decomposition of $G$. Then both the lattices {\rm
Sol($G$)} and {\rm QSol($G$)} are chains of length $\a_k$. More
precisely, under the above notation, we have
$${\rm QSol}(G):\{1\}=\Phi_{\a_k}(G)\subset\Phi_{\a_k-1}(G)\subset ... \subset\Phi_1(G)\subset
\Phi_0(G)=G$$and
$${\rm Sol}(G):\{1\}=\delta(\Phi_0(G))\subset\delta(\Phi_1(G))\subset\delta(\Phi_2(G))\subset
... \subset\delta(\Phi_{\a_k}(G))=G.$$}

Notice that an alternative way of proving Theorem 3.3 can be
inferred from classification of finite abelian groups and the
types of their subgroups. We also observe that for such a group
$G$ the lattice {\rm Sol($G$)} (as well as {\rm QSol($G$)}) is
isomorphic to the lattice $\pi_e(G)$ of element orders of $G$,
because they are direct product of chains of the same length.

\bigskip\noindent{\bf Example.} The lattices
QSol($\mathbb{Z}_2\times\mathbb{Z}_4$) and
Sol($\mathbb{Z}_2\times\mathbb{Z}_4$) associated to the finite
abelian 2-group $\mathbb{Z}_2\times\mathbb{Z}_4$ consist of the
following chains
$${\rm QSol}(\mathbb{Z}_2\times\mathbb{Z}_4):\hspace{3mm}\{1\}\subset\Phi(\mathbb{Z}_2\times\mathbb{Z}_4)\cong\mathbb{Z}_2\subset\mathbb{Z}_2\times\mathbb{Z}_4$$and
$${\rm Sol}(\mathbb{Z}_2\times\mathbb{Z}_4):\hspace{3mm}\{1\}\subset\delta(\Phi(\mathbb{Z}_2\times\mathbb{Z}_4))\cong\mathbb{Z}_2\times\mathbb{Z}_2\subset\mathbb{Z}_2\times\mathbb{Z}_4,$$respectively.

\bigskip\noindent{\bf Remark.} The lattice
QSol($\mathbb{Z}_4\times\mathbb{Z}_4$) is a chain of length 2,
too. So, we have ${\rm QSol}(\mathbb{Z}_2\times\mathbb{Z}_4)\cong
{\rm QSol}(\mathbb{Z}_4\times\mathbb{Z}_4)$. This shows that there
exist non-isomorphic finite groups $G_1$ and $G_2$ such that ${\rm
QSol}(G_1)\cong {\rm QSol}(G_2)$. Moreover, we remark that each of
the following conditions ${\rm QSol}(G_1)\cong {\rm QSol}(G_2)$ or
${\rm Sol}(G_1)\cong {\rm Sol}(G_2)$ is not sufficient to assure
the lattice isomorphism $L(G_1)\cong L(G_2)$.

\bigskip\noindent{\bf Corollary 3.4.} {\it The lattice {\rm QSol($G$)} associated to a finite
abelian group $G$ is a direct product of chains, and therefore it
is decomposable and distributive.}
\bigskip

The properties of the lattice {\rm QSol($G$)} in Corollary 3.4 are
also satisfied by other classes of finite groups $G$. One of them,
which is closely connected to abelian groups, is the class of
hamiltonian groups, that is the finite nonabelian groups all of
whose subgroups are normal. Such a group $H$ can be written as the
direct product of the quaternion group
$$Q_8=\langle x, y \mid x^4=y^4=1, yxy^{-1}=x^{-1}\rangle,$$an elementary
abelian 2-group and a finite abelian group $A$ of odd order, that
is
$$H\cong Q_8 \times \mathbb{Z}_2^n \times A.$$Since $Q_8 \times \mathbb{Z}_2^n$
and $A$ are of coprime orders, the structure of {\rm QSol($H$)}
can be completely described in view of the above results.

\bigskip\noindent{\bf Corollary 3.5.} {\it Let $H\cong Q_8 \times \mathbb{Z}_2^n \times
A$ be a finite hamiltonian group. Then the lattice {\rm QSol($H$)}
is distributive. More precisely, it possesses a direct
decomposition of type
$${\rm QSol}(H)\cong {\rm QSol}(Q_8 \times \mathbb{Z}_2^n) \times {\rm
QSol}(A),$$where ${\rm QSol}(Q_8 \times \mathbb{Z}_2^n)$ is a
chain of length {\rm 3} and ${\rm QSol}(A)$ is a direct product of
chains.}
\bigskip

As show our previous examples, the lattices {\rm Sol($G$)} and
{\rm QSol($G$)} associated to an arbitrary finite group $G$ are
distinct, and this remark remains also valid even for finite
abelian groups. So, the next question is natural: \textit{which
are the finite abelian groups $G$ satisfying {\rm Sol($G$)}={\rm
QSol($G$)}}? By Theorem 3.3, for an abelian $p$-group
$G\cong\prod_{i=1}^k \Z_{p^{\a_i}}$ we have {\rm Sol($G$)}={\rm
QSol($G$)} if and only if the Frattini series and the dual
Frattini series of $G$ coincide, that is $\a_1=\a_2=...=\a_k$.
Clearly, this leads to the following result.

\bigskip\noindent{\bf Corollary 3.6.} {\it Let $G$ be a finite abelian group. Then
{\rm Sol($G$)}{\rm=}{\rm QSol($G$)} if and only if every Sylow
subgroup of $G$ is of type $\prod_{i=1}^k \Z_{p^{\a}}$ for some
positive integer $\a$.}
\bigskip

Finally, Theorem 3.3 can be used to determine the finite groups
that are quotient solitary free. We already know that such a group
$G$ is either perfect or abelian. Note that we were unable to give
a precise description of quotient solitary free perfect groups.
For an abelian group $G$, it is clear that {\rm QSol($G$)} becomes
a chain of length 1 if and only if $G$ is elementary abelian.
Hence the following theorem holds.

\bigskip\noindent{\bf Theorem 3.7.} {\it If a finite group is quotient solitary
free, then it is perfect or elementary abelian. In particular, a
finite nilpotent group is quotient solitary free if and only if it
is elementary abelian.}
\bigskip

\section{Conclusions and further research}

All above results show that the study of some new posets of
subgroups associated to a group $G$, as Sol($G$) and QSol($G$), is
an interesting aspect of subgroup lattice theory. Clearly, it can
be continued by investigating some other properties of these
lattices and can be also extended to more large classes of groups.
This will surely constitute the subject of some further research.
\bigskip

We end this paper by indicating several open problems concerning
these two lattices.
\bigskip

\noindent{\bf Problem 4.1.} Which are the (finite) groups $G$ such
that Sol($G$) and QSol($G$) satisfy a certain lattice-theoretical
property: modularity, distributivity, complementation,
pseudocomplementation, ... and so on?
\bigskip

\noindent{\bf Problem 4.2.} Determine precisely the (finite)
groups $G$ such that QSol($G$)\hspace{1mm}= $N(G)$,
QSol($G$)\hspace{1mm}= Char($G$) and QSol($G$)\hspace{1mm}=
Sol($G$), respectively.
\bigskip

\noindent{\bf Problem 4.3.} Study the containment of some other
important subgroups of a group $G$ (as $\Phi(G)$, $Z(G)$, $F(G)$,
... and so on) to the lattices Sol($G$) and QSol($G$).
\bigskip

\noindent{\bf Problem 4.4.} What can be said about two groups
$G_1$ and $G_2$ for which we have
Sol($G_1$)$\hspace{1mm}\cong\hspace{1mm}$Sol($G_2$) or
QSol($G_1$)$\hspace{1mm}\cong\hspace{1mm}$QSol($G_2$)?
\bigskip

\bigskip\noindent{\bf Acknowledgements.} The author is grateful to the reviewers for
their remarks which improve the previous version of the paper.

\vspace*{5ex}\small

\hfill
\begin{minipage}[t]{5cm}
Marius T\u arn\u auceanu \\
Faculty of  Mathematics \\
``Al.I. Cuza'' University \\
Ia\c si, Romania \\
e-mail: {\tt tarnauc@uaic.ro}
\end{minipage}

\end{document}